\documentclass [a4paper, 12pt] {amsart}

\usepackage {hyperref}

\usepackage [utf8] {inputenc}

\usepackage{mathtext}

\usepackage {version}

\usepackage {amsmath}
\usepackage {amssymb}
\usepackage {amsthm}

\usepackage {graphicx}

\usepackage {appendix}

\usepackage {url}

\usepackage{mathtext}
\usepackage{amsmath}
\usepackage{amsfonts}
\usepackage{amssymb}
\usepackage{mathrsfs}
\usepackage{amsthm}

\usepackage {needspace}

\excludeversion {svntags}

\newtheorem {theorem} {Theorem}

\newtheorem {lemma} [theorem] {Lemma}
\newtheorem {proposition} [theorem] {Proposition}
\newtheorem {definition} [theorem] {Definition}
\newtheorem {corollary} [theorem] {Corollary}

\newcommand {\apclass} [1] {\ensuremath{\mathrm A_{#1}}}

\newcommand {\lclass} [2] {\ensuremath{\mathrm L_{#1} \left( #2 \right) }}
\newcommand {\lsclass} [1] {\ensuremath{\mathit l^{#1} }}

\newcommand {\lsclassl} [1] {\ensuremath{\mathit l_\lambda^{#1} }}

\newcommand {\hclass} [2] {\ensuremath{\mathrm H_{#1} \left( #2 \right) }}
\newcommand {\lclassg} [1] {\ensuremath{\mathrm L_{#1}}}

\newcommand {\hclassg} [1] {\ensuremath{\mathrm H_{#1}}}

\newcommand {\BMO} {\ensuremath {\mathrm {BMO}}}

\newcommand {\AK} {\ensuremath {\mathrm {AK}}}

\newcommand {\nplus} {\ensuremath {\mathrm {N}^+}}

\DeclareMathOperator* {\esssup} {ess\,sup}
\DeclareMathOperator* {\supp} {supp}

\DeclareMathOperator* {\sign} {sign}

\excludeversion {comment}

\newcommand {\weightu} {\ensuremath {\mathit u}}

\newcommand {\weightw} {\ensuremath {\mathit w}}

\def\XXint#1#2#3{{\setbox0=\hbox{$#1{#2#3}{\int}$ }
\vcenter{\hbox{$#2#3$ }}\kern-.6\wd0}}

\begin {document}

\title [BMO-regularity and interpolation of Hardy-type spaces] {On K-closedness,
BMO-regularity
\\
and real interpolation
\\
of Hardy-type spaces}
\author [D.~V.~Rutsky] {Dmitry V. Rutsky}
\email {rutsky@pdmi.ras.ru}
\date {\today}
\address {St.Petersburg Department
of Steklov Mathematical Institute RAS
27, Fontanka
191023 St.Petersburg
Russia}


\maketitle

\begin {abstract}
Let $(X, Y)$ be a suitable couple of quasi-Banach lattices of mea\-surable functions on $\mathbb T \times \Omega$,
and let $(X_A, Y_A)$ be the couple of the corresponding Hardy-type spaces.
We show that $(X_A, Y_A)$ is $\mathrm K$-closed in $(X, Y)$ if and only if $(X, Y)$ is $\BMO$-regular
for a general couple of Banach lattices
having the Fatou property when $\Omega$ is a discrete measurable space.
Furthermore, we establish this equivalence if $X$ is allowed to be quasi-Banach
but $Y$ is assumed to be $p$-convex with some $p > 1$ (here $\Omega$ is arbitrary).
We also show under certain mild restrictions that the ``good interpolation'' formula
$$
\left(X_A, \hclassg {q}\right)_{\theta, p} = \left[\left(X, \lclassg {q}\right)_{\theta, p}\right]_A
$$
holds true if and only if $X$ is $\BMO$-regular.

\end {abstract}

\section {Description of the main results}

\label {mainresults}

The interpolation properties of Hardy-type spaces are of considerable interest; see, e.~g., \cite {kisliakov1999},
\cite [Chapter~7] {kislyakovkruglyak2012}, \cite {thesis}.
For convenience and clarity we begin by briefly stating the results;
a detailed overview of the relevant background with appropriate references is provided in Section~\ref {introduction} below.

We work with quasi-Banach lattices $X$ of measu\-rable functions
on the measurable space $\mathbb T \times \Omega$,
where $(\Omega, \mu)$ is some $\sigma$-finite measurable space
(see Section~\ref {somepreparations} below for the common definitions)
and the associated \emph {Hardy-type spaces}
$$
X_A = \{f \in X \mid f (\cdot, \omega) \in \nplus \text { for a.~e. $\omega \in \Omega$}\},
$$
where $\nplus$ is the boundary Smirnov class.  For example, the Lebesgue spaces $\lclassg {p}$, $0 < p \leqslant \infty$
yield the usual Hardy spaces $\left[\lclassg {p}\right]_A = \hclassg {p}$, but this definition also yields
the Hardy-Lorentz spaces $\hclassg {p, q}$, the weighted Hardy spaces $\hclass {p} {\weightw}$,
the variable exponent Hardy spaces $\hclassg {p (\cdot)}$, the vector-valued Hardy spaces $\hclass {p} {\lsclass {q}}$
and many others.

To avoid degeneration, we usually assume that $X$ satisfies \emph {property $(*)$}:
for any $f \in X$ such that $f \neq 0$ there exists a majorant $g \geqslant |f|$
such that $\|g\|_X \leqslant C \|f\|_X$ and $\log g (\cdot, \omega) \in \lclassg {1}$ for a.~e. $\omega \in \Omega$
with some constant $C$ independent of $f$.
If $X$ also satisfies the Fatou property then $X_A$ is a closed subspace of $X$;
see, e.~g., \cite [\S 1] {kisliakov2002en}.

\begin {definition}
\label {stabilitydef}
Let $(X, Y)$ be a couple of quasi-Banach lattices of measurable functions on $\mathbb T \times \Omega$.
We say that the couple $(X_A, Y_A)$ has stable interpolation with respect to an interpolation functor
$\mathcal F$ if it satisfies the ``good interpolation formula''
\begin {equation}
\label {gif}
\mathcal F \left(\left(X_A, Y_A\right)\right) = \left[\mathcal F \left((X, Y)\right) \right]_A.
\end {equation}
\end {definition}

An important case of stability for the real interpolation is outlined by the following property.
\begin {definition}
\label {akstdef}
A couple $(X, Y)$ is called $\AK$-stable
if the couple $(X_A, Y_A)$ is $\mathrm K$-closed in $(X, Y)$, that is,
if for any $f \in X$ and $g \in Y$ such that $f + g \in X_A + Y_A$
there exist some $F \in X_A$ and $G \in Y_A$ satisfying
$f + g = F + G$, $\|F\|_X \leqslant c \|f\|_X$ and $\|G\|_Y \leqslant c \|g\|_Y$
with some constant $c$ independent of $f$ and $g$.
\end {definition}
Informally, $\AK$-stability means that we may replace arbitrary measu\-rable decompositions
of functions
from $X_A + Y_A$ by some analytic decompositions with a good control on the norm.
As one would imagine, this property has many applications.
In particular, it is easy to see that \eqref {gif} holds true for all $\AK$-stable couples $(X, Y)$
and real interpolation functors~$\mathcal F$.

In the interesting cases
$\AK$-stability is difficult to verify directly; however, as we will see in a moment, it is closely related
to the following property.
\begin {definition}
\label {bmordef}
A lattice $X$ is said to be BMO-regular
if for any nonzero $f \in X$ there exists a majorant $u \geqslant |f|$ such that
$\|u\|_X \leqslant m \|f\|_X$ and $\left\|\log u (\cdot, \omega)\right\|_{\BMO} \leqslant C$ for
a.~e. $\omega \in \Omega$ with some constants $m$ and $C$ independent of $f$.
\end {definition}
\begin {definition}
\label {bmorcdef}
A couple
$(X, Y)$ is called \emph {BMO-regular} if for all nonzero $f \in X$ and $g \in Y$
there exist some majorants
$u \geqslant |f|$ and $v \geqslant |g|$ such that $\|u\|_X \leqslant m \|f\|_X$, $\|v\|_Y \leqslant m \|g\|_Y$
and $\left\|\log \frac {u (\cdot, \omega)} {v (\cdot, \omega)}\right\|_{\BMO} \leqslant C$ for a.~e. $\omega \in \Omega$
with some constants $C$ and $m$ independent of $f$ and $g$.
\end {definition}
This property describes a class of couples of lattices that are nice in a certain natural sense
pertaining to harmonic analysis,
and as such $\BMO$-regularity is fairly well understood; see, e.~g., \cite {rutsky2011en}.
In many particular cases it is possible to give simple characterizations for $\BMO$-regular couples.
For instance, 
a weighted Lebesgue space $\lclass {p} {\weightw}$
with the norm $\|f\|_{\lclass {p} {\weightw}} = \left\|f \weightw^{-1}\right\|_{\lclassg {p}}$ is $\BMO$-regular
if and only if $\log \weightw (\cdot, \omega) \in \BMO$ uniformly in a.~e. $\omega \in \Omega$,
and a couple $\left(\lclass {p_0} {\weightw_0}, \lclass {p_1} {\weightw_1}\right)$ is $\BMO$-regular if and only if
$\log \frac {\weightw_0 (\cdot, \omega)} {\weightw_1 (\cdot, \omega)} \in \BMO$ 
uniformly in a.~e. $\omega \in \Omega$.

At first glance it may seem that $\BMO$-regularity has nothing to do with $\AK$-stability.
However, $\BMO$-regularity implies $\AK$-stability by~\cite [Theorem~3.3] {kisliakov1999}.
It has long been suspected that $\BMO$-regularity is not only sufficient but is also necessary for $\AK$-stability,
but the equivalence 
under a natural generality
remained elusive;
see Section~\ref {introduction} below for a detailed background. 
Now we are finally able to provide such a result.
\begin {theorem}
\label {akbakeq}
Let $X$ and $Y$ be Banach lattices of measurable functions on
$\mathbb T \times \Omega$ satisfying the Fatou property and property $(*)$.
Suppose also that $\Omega$ is a discrete space.
Then $(X, Y)$ is $\AK$-stable if and only if $(X, Y)$ is $\BMO$-regular.
\end {theorem}
The proof of Theorem~\ref {akbakeq} 
is given in Section~\ref {akbakeqs} below;
we establish that under the stated conditions
$\AK$-stability implies the
bounded $\AK$-stability property, which is known to be equivalent to $\BMO$-regularity
by \cite [Theorem~4] {rutsky2013t}.
The idea and many of the details of the argument are essential\-ly identical to the main result
of \cite {kisliakovrutsky2012en}: we use a fixed point theorem in order to derive the existence of a decomposition
of the required form from a weaker property.
Although this approach has been known to the author at the time the question about the equivalence
between $\AK$-stability and bounded $\AK$-stability was raised in~\cite {rutsky2010en},
the technique developed by that time proved to be inadequate.
The Fan--Kakutani fixed point theorem, which was successfully applied before to similar problems
(see \cite {kisliakov2002en}, \cite {rutsky2010en}, \cite {rutsky2011en}, \cite {kisliakovrutsky2012en}),
seems to be inapplicable to the task at hand, and we have to take advantage of a much more potent
Powers's fixed point theorem~\cite {powers1972}.

The assumption that the space $\Omega$ is discrete is crucial for the particular fixed point argument
in the proof
of Theorem~\ref {akbakeq}, as is the assumption that both $X$ and $Y$ are Banach lattices,
and it is presently unclear how these restrictions might be lifted.
The discreteness assumption in particular, although apparently not terribly restrictive,
still seems to be very odd and unfortunate by contrast with the rest of the theory working fine for arbitrary $\Omega$,
including the rest of the results of the present work.

If one of the lattices
has nontrivial convexity,
we may take any quasi-Banach lattice for the other and still establish the equivalence.
\begin {theorem}
\label {akbakeqpc}
Let $X$ be a quasi-Banach lattice and $Y$ be a $p$-convex Banach lattice with some $p > 1$,
both spaces being lattices of measurable functions on
$\mathbb T \times \Omega$ satisfying the Fatou property and property $(*)$.
Then $(X, Y)$ is $\AK$-stable if and only if $(X, Y)$ is $\BMO$-regular.
\end {theorem}
The proof of Theorem~\ref {akbakeqpc} is given in Section~\ref {akbakeqpcs} below.
This result follows from
\cite [Theorem~2] {rutsky2013t} and is thus independent of Theorem~\ref {akbakeq}.
The main idea is that we can replace $X$ by a real interpolation space $(X, Y)_{\theta, p}$, which has
the required convexity if $\theta$ is sufficiently close to $1$;
we arrive at the conditions of \cite [Theorem~2] {rutsky2013t}
after a reasonably short chain of reductions with the help of lattice multiplication and duality.
\begin {corollary}
\label {akbakeqpco}
Let $(X, Y)$ be a couple of Banach lattices of measurable functions on
$\mathbb T \times \Omega$ satisfying the Fatou property and property $(*)$.
Suppose that $Y$ is $q$-concave with some $1 \leqslant q < \infty$.
Then $(X, Y)$ is $\AK$-stable if and only if $(X, Y)$ is $\BMO$-regular.
\end {corollary}
This immediately follows from Theorem~\ref {akbakeqpc}, since $Y'$ is a $q'$-convex lattice and
couples $(X, Y)$ and $(X', Y')$ are $\AK$-stable and $\BMO$-regular only simultaneously by
\cite [Lemma~7] {kisliakov2002en} and \cite [Theorem~5.8] {rutsky2011en}.
\begin {corollary}
\label {akbmoeql}
Let $X$ be a Banach lattice of measurable functions on $\mathbb T \times \Omega$
satisfying the Fatou property and property $(*)$,
and let $1 \leqslant p \leqslant \infty$.
The couple $(X, {\rm L}_{p})$ is $\AK$-stable if and only if $X$ is $\BMO$-regular.
If $p > 1$ then the same is true for all quasi-Banach lattices $X$.
\end {corollary}
This corollary was established in~\cite [Theorem~3] {rutsky2013t} for Banach lattices~$X$ and
under an awkward assumption that
$p \in \{1, 2, \infty\}$, which we thus lift.
The case $p > 1$ is a direct consequence of Theorem~\ref {akbakeqpc},
and the case $p = 1$ follows from Corollary~\ref {akbakeqpco}.

In fact, for most couples $(X, \lclassg {q})$ it is possible to obtain a stronger result,
namely the equivalence of
the stability \eqref {gif} for the usual real interpolation functors
and $\BMO$-regularity.
\begin {theorem}
\label {bmogrilp}
Let $X$ be a Banach lattice of measurable functions on $\mathbb T \times \Omega$
satisfying the Fatou property and property $(*)$,
and let $1 \leqslant p, q < \infty$.
Suppose that both $X$ and $X'$ have order continuous norm,
and either $q > 1$, or $X$ is $r$-convex with some $r > 1$.
Then
\begin {equation}
\label {rigip}
\left(X_A, \hclassg {q}\right)_{\theta, p} = \left[\left(X, \lclassg {q}\right)_{\theta, p}\right]_A
\end {equation}
for some (equivalently, for all) $0 < \theta < 1$ if and only if $X$ is $\BMO$-regular.
\end {theorem}
The proof of Theorem~\ref {bmogrilp} is given in Section~\ref {bmogrilps} below;
it uses the complete form
of the seminal result \cite [Theorem~5.12] {kalton1994} of N.~Kalton (see Section~\ref {introduction} below).
The super\-reflexivity assumptions imposed in \cite {kalton1994} are very similar to that of Theorem~\ref {bmogrilp}.
The case $Y = \lclassg {q}$ lends itself to a rather simple treatment:
as in the main result of \cite {rutsky2013t}, after certain reductions we can
take advantage of the formula $(Z, Z')_{\frac 1 2, 2} = \lclassg {2}$
for a suitable lattice $Z$.
It is yet to be seen whether such an equivalence holds true for general couples $(X, Y)$.

\section {Historical remarks}

\label {introduction}

By necessity we only sketch out some highlights and allow (hopefully minor) simplifications and omissions;
see also
\cite {kisliakov1999},
\cite [Chapter~7] {kislyakovkruglyak2012}, \cite {thesis}.

For the classical couples $\left(\hclassg {p}, \hclassg {q}\right)$ with $1 \leqslant p, q \leqslant \infty$
(where $p = 1$ and $q = \infty$ are the interesting cases, since the corresponding Hardy spaces are not complemented in
$\lclassg {1}$ and~$\lclassg {\infty}$)
P.~Jones in \cite {jones1983} established the stability of interpolation~\eqref {gif} for arbitrary interpolation
functors $\mathcal F$;
see also \cite [\S 5.10] {bennetsharpley}.
The proof uses certain constructive solutions for $\bar\partial$-equations with Carleson measure data
and is thus rather involved.

Apparently, for the first time the importance of the $\AK$-stability\footnote {The term \emph {AK-stability},
a convenient
shorthand for ``analytic K-stability'', was introduced
in \cite {kisliakov2003}.}
property was clearly
articulated by G.~Pisier in \cite {pisier1992}, \cite {pisier1992b}, \cite {pisier1993}, who, in particular,
proved in a rather elementary way
the $\AK$-stability of $\left(\lclassg {p}, \lclassg {\infty}\right)$ for all $p > 0$ and
of $\left(\lclass {p_0} {\lsclass {q_0}}, \lclass {p_1} {\lsclass {q_1}}\right)$ for $1 \leqslant p_j, q_j \leqslant \infty$,
and also gave a simple proof of the Grothendieck theorem for the disc algebra with the help of the latter.

After a few further developments (that we will mention shortly)
came a landmark contribution by N.~Kalton \cite {kalton1994}, who did a thorough study of the stability of
the complex interpolation in a rather general setting.
Apparently for the first time, he introduced\footnote { To be precise, for the first time
Definition~\ref {bmorcdef} in this particular form was given in \cite {kisliakov1999};
in \cite {kalton1994}, a notion of \emph {BMO-direction} was introduced
in line with a rather complicated construction
(in pursuit of generality) that can be shown (under the usual assumptions) to be equivalent to
the $\BMO$-regularity in the light of the later developments.}
the notion of $\BMO$-regularity and
pointed out most of its essential properties.
The $\BMO$-regularity theory was further extended and refined in
\cite {kisliakov1999}, \cite {kisliakov2002en}, \cite {kisliakov2003}, \cite {rutsky2010en}, \cite {rutsky2011en}.


In a series of examples preceding \cite {kalton1994} it was gradually discovered
that $\BMO$ has something to do with the interpolation of Hardy-type spaces.
M.~Cwikel, J.~McCarthy and T.~Wolff showed in~\cite {cwikel1ibw} that\footnote {
For simplicity, we gloss over the fact that some of the results mentioned
below require that the weights $\weightw$ that appear in them satisfy $\log \weightw \in \lclassg {1}$,
in accordance with property $(*)$ usually assumed in the results for the general lattices.}
$\hclass {p} {\weightw_0^{1 - \theta} \weightw_1^\theta}$ is an interpolation space of type $\theta$, $0 < \theta < 1$,
for the couple
$\left(\hclass {p} {\weightw_0}, \hclass {p} {\weightw_1}\right)$
if and only if $\log \frac {\weightw_0} {\weightw_1} \in \BMO$ using the corona theorem.
Later S.~Kisliakov and Q.~Xu established in \cite {kisliakovxu1994} that
the couple $\left(\hclass {p_0} {\weightw_0}, \hclass {p_0} {\weightw_1}\right)$ is $\AK$-stable if and only if
$\log \frac {\weightw_0} {\weightw_1} \in \BMO$ and made some vector-valued generalizations formulated in terms
of a majorization property, which was subsequently seen (in \cite {kisliakov1998}) to be equivalent to $\BMO$-regularity.

In~\cite {kalton1994} N.~Kalton proved that under certain superreflexivity assumptions
a couple $(X, Y)$ satisfies \eqref {gif} for the complex interpolation functor
$\mathcal F = (\cdot, \cdot)_\theta$ with some $0 < \theta < 1$ if and only if $(X, Y)$ is $\BMO$-regular,
thus completely characterizing for the first time an interpolation stability property for Hardy-type spaces
in a rather general setting.
The proof of the ``if'' part,
given by \cite [Theorem~5.7] {kalton1994}, amounts
to a very simple explicit construction\footnote {This result was later generalized in~\cite [Theorem~3.4] {kisliakov1999};
see also~\cite [Corollary~2] {kisliakovxu2000}.}
(which can easily be made to work
for $\BMO$-regular couples whenever both lattices have order continuous norm).
The ``only if'' part, on the other hand, is rather involved; together with the early version of the $\BMO$-regularity
theory, it is based on Kalton's remarkable tools
for computations with K\"othe spaces, which have attracted significant attention over the years
(see, e.~g., \cite {cwikelmilmanrochberg2014}), and on some deep insights into the stability of complex interpolation.

Major progress was achieved in a review \cite {kisliakov1999} by S.~Kisliakov, who demonstrated, among other things,
that every $\BMO$-regular couple is $\AK$-stable, and provided a
``method solving all interpolation problems'' for couples of weighed Hardy spaces, proving in a rather elementary
manner that every $\BMO$-regular couple of weighted Lebesgue spaces satisfies~\eqref {gif} and thus resolving just
as much in this setting as P.~Jones did some 15 years earlier for the classical case; see Theorem~\ref {goodinterphp} in
Section~\ref {concrem} below. 

The question about the relationship between $\AK$-stability and $\BMO$-regularity was (somewhat implicitly)
raised in \cite {kisliakov1999}.
It was well known at that time that this is the case for couples of weighted Lebesgue spaces,
but Kalton's result suggested that this might also be true in a much more general setting.

Indeed, in~\cite {kisliakov2002en} it was shown that $\BMO$-regularity is necessary for $\AK$-stability
for couples of a special form $\left(X (\lsclass {p}), \lclass {\infty} {\lsclassl {\infty}}\right)$
with an additional variable, but without any significant restrictions on the lattice~$X$.
The proof introduced several new ideas; the arguments making use of an additional variable
with a power weight and the application of a fixed point theorem are of particular importance.
Fixed point arguments later
proved to be instrumental in establishing certain important properties of $\BMO$-regularity
(see~\cite {rutsky2010en}, \cite {rutsky2011en}).

The equivalence established in~\cite {kisliakov2002en}
was later extended to some other similar cases in~\cite {rutsky2010en}, where also
a \emph {bounded AK-stability} property, which is a stronger and in many respects more convenient
modification of the $\AK$-stability property,
was explicitly introduced and studied to some extent.
In \cite {thesis}, it was shown (again with the help of the bounded $\AK$-stability)
that the equivalence of $\AK$-stability and $\BMO$-regularity
takes place for couples $\left(\lclassg {p (\cdot)}, \lclassg {\infty}\right)$
with a piecewise smooth exponent~$p (\cdot)$.  Although in retrospect
these results
do not seem to be particularly impressive and were not directly useful in the end,
they still helped to build up some confidence and refine certain techniques.

In \cite {kisliakovrutsky2012en},
a technique was developed for a question related to the corona theorem.
In the present work essentially the same technique is used to establish the equivalence between $\AK$-stability and
the bounded $\AK$-stability in the proof of Theorem~\ref {akbakeq}.

More recently in \cite {rutsky2013t}, new insights made it possible to prove the equivalence
of $\AK$-stability and $\BMO$-regularity for general couples of 2-convex lattices $(X, Y)$ and in some similar cases.
Firstly, it was realized that $\AK$-stability of a couple of the form $(Z, Z')$
implies stability for complex interpolation of a related couple
$\left(V_{\theta_1}, V_{\theta_2}\right)$ of spaces $V_{\theta} = (Z, \lclassg {2})_{\theta, 2}$,
which allows one to take advantage of Kalton's results and establish $\BMO$-regularity of $V_\theta$.
The proof of Theorem~\ref {bmogrilp} in Section~\ref {bmogrilps} uses a similar reduction.
Secondly, it was shown that $\BMO$-regularity of $V_\theta$ implies $\BMO$-regularity of $Z$,
in a complicated argument that introduced an additional variable and exploited a symmetry to infer
$\BMO$-regularity of 
$Y^{\frac 1 2}$
from that of
$(Y, \lclassg {\infty})_{\frac 1 2, \infty}$.
In the computations needed to make the necessary reduction,
a three-lattice reiteration formula was used.
The validity of such formulae is based on certain nontrivial facts about the Sparr interpolation spaces;
similar computations appear in the proof of Proposition~\ref {genridiv} below.

\section {Preparations}

\label {somepreparations}

We begin by stating a very general fixed point theorem from~\cite {park1998} used in the main result;
for a good general reference on the fixed point theory see, e.~g., \cite {granasdugundji2003}.
Luckily, this impressive result is not hard to get hold of for our purposes,
even though the relevant theory is rather complicated and its key elements may not be familiar or even quite readily
explainable to an interested reader coming from analysis.

Suppose that $X$ and $Y$ are topological spaces.  A set-valued map 
$T : X \to 2^Y$ is called \emph {closed} if its graph is closed in $X \times Y$.
$T$ is said to be \emph {upper semicontinuous} if for any
closed set $B \subset Y$
its preimage
$$
T^{-1} (B) = \left\{x \in X \mid T (x) \cap B \neq \emptyset \right\}
$$
is also closed.
There is a more natural equivalent definition: $T$ is upper semicontinuous if and only if
for any open set $U \subset Y$ the set
$$
\{x \in X \mid T (x) \subset U\}
$$
is also open.
Thus a composition of of upper semicontinuous maps is also upper semicontinuous.
It is easy to see that if $Y$ is a regular topological
space\footnote {That is, we can separate a point from a closed set not containing it by a couple of open neighbourhoods;
it is well known that any Hausdorff topological vector space is regular.}
and the values of $T$ are closed then $T$ is upper semicontinuous if and only if $T$ is a closed map.
$T$ is called \emph {compact} if the closure of its image $\overline {T (X)}$ is compact in $Y$.
Observe that a composition of compact maps (and even a composition of a compact map with any map)
defined on a Hausdorff compact set is also compact.
For the notion and the definition of an \emph {acyclic} topological space we (by necessity)
refer the reader to \cite {granasdugundji2003}
and to the various algebraic topology textbooks;
in the present work we will only use the simple fact
that convex sets of a topological vector space are acyclic.
$T$ is called an \emph {acyclic map} if $T$ is upper semicontinuous and its values are compact and acyclic.

A nonempty set $X \subset E$ in a linear topological space $E$ is called \emph {admissible} (in the sense of Klee)
if for any compact set $K \subset X$ and any open set $V \subset E$, $0 \in V$, there exists
a continuous map $h : K \to X$ such that
$x - h (x) \in V$ for all $x \in K$ and $h (K)$ is contained in a finite-dimensional subspace $L \subset E$.
In other words, $X$ is admissible if any compact set $K \subset X$ can be continuously and uniformly approximated
by a family of finite-dimensional sets of $X$.
In particular, any non\-empty convex set of a locally convex linear topological space is admissible.

Let $X$ be a nonempty convex set in a linear topological space $E$,
and let $Y$ be another linear topological space.
A set $P \subset X$ is called a \emph {polytope} if $P$ is the convex hull of a finite set in $X$.
A map $F : X \to 2^Y$ belongs to the \emph {``better'' admissible class} $\mathfrak B (X, Y)$ if and only if
for any polytope $P \subset X$ and any continuous function $f : F (P) \to P$ the composition
$f \circ F|_P : P \to P$ has a fixed point.
Observe that \emph {admissibility} refers in this notion to the existence of fixed points in a restricted sense.
The class $\mathfrak B (X, Y)$ encompasses a large number of particular classes of maps that are known to have fixed points;
in the present work we will only use the fact that this class contains finite compositions of acyclic maps.
The corresponding fixed point theorem was established in~\cite {powers1972} (see also \cite [\S 19.9] {granasdugundji2003}
for the statement in context), and 
it is possible to use it directly with minor adaptations; the powerful result~\cite {park1998}, however, allows us to keep
the necessary topological explanations to a minimum.
\begin {theorem} [{\cite [Corollary~1.1] {park1998}}]
\label {parkfptc}
Let $E$ be a Hausdorff topological vector space, and let $X \subset E$ be an admissible convex set.
Then any closed compact map $F \in \mathfrak B (X, X)$ has a fixed point.
\end {theorem}

Let $(\Omega, \mu)$ be a $\sigma$-finite measurable space, and let $m$ be
the normed Lebesgue measure on the unit circle $\mathbb T$.
A quasi-normed lattice of mea\-surable functions $X$ on $\mathbb T \times \Omega$
is a quasi-normed space of measurable functions $X$
in which the norm is compatible with the natural order; that is, if $|f| \leqslant g$ a.~e. for some function $g \in X$
then $f \in X$ and $\|f\|_X \leqslant \|g\|_X$.
For simplicity we only work with lattices $X$ such that $\supp X = \mathbb T \times \Omega$ up to a set of measure $0$.
For more detail on the normed lattices and their properties see, e.~g., \cite [Chapter~10] {kantorovichold}.

A function $f$ on $X$ is said to be \emph {order continuous} if for any sequence $x_n \in X$
such that $\sup_n |x_n| \in X$ and $x_n \to 0$ a.~e. one also has $f (x_n) \to 0$.
If $X$ is Banach then any order continuous functional $f$ on $X$
has an integral representation
$f (x) = \int x y_f$ for some measurable function $y_f$ which can be identified with $f$.
The set of all such functionals $X'$ is also a Banach lattice with
a norm defined by $\|f\|_{X'} = \sup_{g \in X, \|g\|_X = 1} \int |f g|$.
Lattice $X'$ is called the \emph {order dual} of the lattice $X$.

A lattice $X$ has the Fatou property if for any
$f_n, f \in X$ such that $\|f_n\|_X \leqslant 1$ and the sequence $f_n$ 
converges to $f$ a.~e. it is also true that $f \in X$ and $\|f\|_X \leqslant 1$.
The Fatou property of a lattice $X$ is equivalent to $(m \times \mu)$-closedness of the unit ball
$B_X$ of the lattice $X$
(here and elsewhere $(m \times \mu)$-convergence denotes the
convergence in measure in any measurable set $E$ such that $(m \times \mu) (E) < \infty$).
If $X$ has the Fatou property then the unit ball of $X_A$ is also closed with respect to the convergence in measure
by the Khinchin-Ostrovski theorem
(\cite [Chapter~II, \S 7.1, Theorem~5] {privaloven}); see \cite [\S 1] {kisliakov2002en}.

By definition, for any $f \in X_A$ and almost all $\omega \in \Omega$
function $f (\cdot, \omega)$ on $\mathbb T$ represents the boundary values of an analytic function on $\mathbb D$;
thus any such $f$ is naturally identified with a function on $\mathbb D \times \Omega$ analytic in the first variable
for almost all fixed values of the second variable.
In the case of a discrete measure $\mu$
we can exploit the topology of uniform convergence on the compact sets of $\mathbb D$
in order to establish closedness of maps;
such an argument appeared before in~\cite [\S 3] {kisliakov2002en}.
\begin {proposition}
\label {latballic}
Let $X$ be a Banach lattice of measurable functions on $\left(\mathbb T \times \Omega, m \times \mu\right)$, where
$\mu$ is a discrete measure consisting of a finite number of point masses.
Suppose that $X$ satisfies the Fatou property and property $(*)$.
Then the closed unit ball $B$ of $X_A$ is compact in the topology\footnote {
One may also take for $\tau$ the equivalent topology induced from the $*$-weak topology of $\hclassg {1}$ by
the map $\Phi : X_A \to \hclassg {1}$, $\Phi (f) = W f$, where $W$ is the outer function defined in the proof;
however, the necessary topological explanations seem to complicate this approach and make it less straightforward.
}
$\tau$ of the uniform convergence
on all compact sets of $\mathbb D \times \Omega$.
\end {proposition}
Since $\tau$ is metrizable, it is sufficient to verify that for any sequence
$f_n \in B$ there is a subsequence converging to some $f \in B$ in $\tau$.
Observe that there exists some $g \in X'$ such that $\|g\|_{X'}  = 1$ and $g > 0$ a.~e.
(see, e.~g., \cite [Proposition~9] {rutsky2010en}).
Lattice $X'$ also satisfies property $(*)$ (see \cite [Lemma~2] {kisliakov2002en}),
so there exists some
$\weightw \in X'$ such that $\weightw \in X'$, $\weightw > g > 0$ a.~e. and
$\log \weightw (\cdot, \omega) \in \lclassg {1}$
for a.~e. $\omega \in \Omega$.  We may assume that $\|\weightw\|_{X'} = 1$.
Thus we can construct an outer function
$W = \exp \left(\log \weightw + i H [\log \weightw]\right)$ such that $|W| = \weightw$ a.~e.; here
$H$ denotes the Hilbert transform acting in the first variable.
It is easy to see that $f_n W$ belongs to the unit ball of the space
$\hclass {1} {\mathbb T \times \Omega}$.  Since measure $\mu$ consists of a finite number of point masses,
the latter space is dual to
$\mathrm C \left(\mathbb T \times \Omega\right) \slash \mathrm C_A \left(\mathbb T \times \Omega\right)$,
so there exists a subsequence $f_{n'} W$ converging in the $*$-weak topology to some
$h \in \hclass {1} {\mathbb T \times \Omega}$,
and therefore $f_{n'} W \to h$ in $\tau$.
Finally, we need to verify that
$f = W^{-1} h \in B$.
Indeed, by a well-known corollary to the Fatou property
(see, e.~g., \cite [Proposition~10] {rutsky2010en} or~\cite [Proposition~3.3] {rutsky2011en})
there exists a sequence
$\varphi_j$ of finite convex combinations of $\{f_{n'}\}_{n' > j}$
such that $\varphi_j \to \varphi$ a.~e. on $\mathbb T \times \Omega$
for some $\varphi \in B$, and we surely have $\varphi_j \to f$ on $\mathbb D \times \Omega$.
By the Khinchin-Ostrovski theorem
the sequence $\varphi_j$ converges to $\varphi$ also in $\tau$.
Thus indeed $f = \varphi \in B$.

\begin {lemma}
\label {logl1isl2}
Let $f \in \lclassg {1}$ and $f \geqslant 1$ almost everywhere.
Then $\log f \in \lclassg {2}$ and  $\|\log f\|_{\lclassg {2}} \leqslant 2 \|f\|_{\lclassg {1}}^{\frac 1 2}$.
\end {lemma}
We only need to observe that $\log f^{\frac 1 2} \leqslant f^{\frac 1 2}$, and so
$$
\int (\log f)^2 = 4 \int \left(\log f^{\frac 1 2}\right)^2 \leqslant 4 \int f.
$$

\begin {proposition}
\label {luscompact}
Let $X$ be a Banach lattice of measurable functions on $\left(\mathbb T \times \Omega, m \times \mu\right)$
with measure $\mu$ consisting of a finite number of point masses, and let $f \in X$.
Suppose that $X$ satisfies the Fatou property and property $(*)$, $f > 0$ almost everywhere and
$\log f (\cdot, \omega) \in \lclassg {1}$ for a.~e. $\omega \in \Omega$.
Then for all $A > 0$ sets
\begin {equation}
\label {vfaset}
V_{X, f, A} = \{\log g \in \lclassg {1} \mid g \geqslant f, \|g\|_X \leqslant A\}
\end {equation}
are compact in the weak topology of $\lclassg {1}$.
\end {proposition}
Indeed, by making a change of measure if necessary we may assume that $m (\mathbb T) = \mu (\Omega) = 1$.
It is well known that the positive part of the unit ball of a Banach lattice is logarithmically convex,
so $V_{X, f, A}$ is a convex set.
Observe first that by the Fatou property $V_{X, f, A}$ is closed with respect to the convergence in measure,
so $V_{X, f, A}$ is also closed in $\lclassg {1}$ and thus weakly closed.
Now, by the Dunford-Pettis theorem it suffices to prove that the set $V_{X, f, A}$
is bounded and uniformly absolutely continuous.
Let $B$ be a measurable set of $\mathbb T \times \Omega$.
We define $\weightw \in X'$ as in the proof of Proposition~\ref {latballic}.
Then $\{\weightw g \mid \log g \in V_{X, f, A}\}$ is a bounded set in $\lclassg {1}$,
and by Lemma~\ref {logl1isl2} we have
$\|\log^+ [\weightw g]\|_{\lclassg {2}} \leqslant 2 A^{\frac 1 2}$ for all $\log g \in V_{X, f, A}$.
Thus
$$
\int_B \log^+ [\weightw g] \leqslant \left\|\log^+ [\weightw g]\right\|_{\lclassg {2}} [(m \times \mu) (B)]^{\frac 1 2} \to 0
$$
as $(m \times \mu) (B) \to 0$, and this convergence is uniform in $\log g \in V_{X, f, A}$.
On the other hand, we also have
$$
\int_B \log^- [\weightw g] \geqslant \int_B \log^- [\weightw f] \to 0
$$
as $(m \times \mu) (B) \to 0$ uniformly in $\log g \in V_{X, f, A}$.
Therefore the set
$$
\{\log [\weightw g] \mid \log g \in V_{X, f, A}\} = \log \weightw + V_{X, f, A}
$$
is bounded and uniformly absolutely continuous, which implies its relative
com\-pac\-tness in the weak topology of $\lclassg {1}$.
It follows that $V_{X, f, A}$ is also compact.

A couple $(X, Y)$ is called \emph {strongly $\AK$-stable} if
for any $f \in X$, $g \in Y$, such that $f + g \in (X + Y)_A$, there exist some
$F \in X_A$ and $G \in Y_A$ such that $f + g = F + G$, $\|F\|_X \leqslant c \|f\|_X$ and $\|G\|_Y \leqslant c \|g\|_Y$
with a constant $c$ independent of $f$ and $g$.
The properties of $\AK$-stability and strong $\AK$-stability are equivalent for couples of Banch lattices satisfying
the Fatou property; see \cite [Lemma~3] {kisliakov2003}.

\section {Proof of Theorem \ref {akbakeq}}

\label {akbakeqs}

Suppose that under the conditions of Theorem~\ref {akbakeqs}
a couple $(X, Y)$ is $\AK$-stable.
By~\cite [Theorem~4] {rutsky2013t} it is sufficient to prove that $(X, Y)$ is \emph {boundedly $\AK$-stable}:
we are given some $f \in X$ and $g \in Y$ and
we need to find some $U \in \hclassg {\infty}$ such that $\|U\|_{\lclassg {\infty}} \leqslant c$ and
$\|g U \|_X \leqslant c \|f\|_X$,
$\|f (1 - U) \|_Y \leqslant c \|g\|_Y$ with a suitable constant $c$ depending only on the $\AK$-stability constant
of the couple $(X, Y)$.
To do that, we construct a certain set-valued map such that its fixed points yield a decomposition of the
required form.

First, we may assume that the measure $\mu$ has only a finite number of point masses.
Indeed, suppose that the result holds true in this case. Let $\alpha \subset \Omega$ be a measurable set
containing only a finite number of the point masses of $\mu$.
Then the sets
$$
U_\alpha = \left\{
U \in \hclassg {\infty} \mid
\|\chi_{\mathbb T \times \alpha} [g U] \|_X \leqslant c \|f\|_X, \|\chi_{\mathbb T \times \alpha} [f (1 - U)] \|_Y \leqslant c \|g\|_Y
\right\}
$$
are nonempty.  It is easy to see that the sets $U_\alpha \subset \lclassg {\infty}$
are convex, bounded and closed with respect to the convergence in measure,
and $U_{\beta} \subset U_{\alpha}$ if $\alpha \subset \beta$.
We may take a nondecreasing sequence $\alpha_j$ such that $\bigcup_j \alpha_j = \Omega$ and apply
\cite [Theorem~3, \S X.5] {kantorovichold} to establish that the set $U_\Omega = \bigcap_j U_{\alpha_j}$ is nonempty,
which implies the existence of a function $U \in U_\Omega$ with the necessary properties.

Excluding some trivial cases allows us to assume that
$f \neq 0$ in $X$ and $g \neq 0$ in $Y$.
By making use of the property $(*)$
we may assume that
$$
\log f (\cdot, \omega), \log g (\cdot, \omega) \in \lclassg {1}
$$
for almost all $\omega \in \Omega$.
Furthermore, homogeneity allows us to assume that $\|f\|_X = 1$.  Suppose that $\lambda = \|g\|_Y$.
We denote by $X + \lambda Y$ the space of the measurable functions
$h$ on $\mathbb T \times \Omega$
having a finite $K$-functional norm
$$
\|h\|_{X + \lambda Y} = K \left(\lambda^{-1}, h; X, Y\right) = \inf_{h = a + b} \left( \|a\|_X + \lambda^{-1} \|b\|_Y \right).
$$
It is easy to see that $X + \lambda Y$ is a Banach lattice satisfying the Fatou property.
We denote by $B_Z$ the closed unit ball of $Z$ for any Banach space $Z$.
Surely
$B_X + \lambda B_Y \subset 2 B_{X + \lambda Y}$.

Let $V = V_{f \vee g, X + \lambda Y, 4}$ be the set  defined by~\eqref {vfaset}; we endow $V$ with the weak
topology of $\lclassg {1}$.
By Proposition~\ref {luscompact} $V$ is weakly compact, and so
$V$ is metrizable since $\lclassg {1}$ is a separable space.
We endow spaces $B_{X_A}$ and $B_{Y_A}$ by the respective topologies of uniform convergence 
on compact sets in $\mathbb D \times \Omega$ and define a map
$\Phi_1 : B_{X_A} \times \lambda B_{Y_A} \to 2^V$ by
\begin {equation*}
\Phi_1 \left((u, v)\right) = \left\{
\log \weightw \in \lclassg {1}
\mid
\weightw \geqslant |u| \vee |v| \vee f \vee g,
\|\weightw\|_{X + \lambda Y} \leqslant 4
\right\}.
\end {equation*}
It is easy to see that the values of $\Phi_1$ are nonempty.
By making use of the logarithmic convexity of $B_{X + \lambda Y}$
it is not difficult to see that the values of $\Phi_1$ are convex.
Let us verify that $\Phi_1$ is upper semicontinuous.\footnote {
We note in passing that since the logarithm is a concave function,
in order to have upper semicontinuity of $\Phi_1$ it is crucial that
it is defined on a set of analytic functions and the measure $\mu$ consists of a finite number of point masses.
For a simple illustrative example consider a sequence $u_n (t) = v_n (t) = 1 + (e^2 - 1) \sign^+ \sin (2 \pi n t)$.
It is easy to see that $\log v_n \to 1$ in the weak topology of $\lclassg {1}$,
but in many cases we have $u_n \to \frac {e^2} 2 > 1$; for example,
in the weak topology of the lattice $\lclassg {2}$.
}
It is sufficient to prove that the graph of $\Phi_1$ is closed.
Indeed, suppose that we are given some sequences
$$
(u_n, v_n) \in B_{X_A} \times \lambda B_{Y_A}
$$
and
$\log \weightw_n \in \Phi_1 \left((u_n, v_n)\right)$
such that $u_n \to u$, $v_n \to v$ and $\log \weightw_n \to \log \weightw$
in the respective spaces; we are to verify that $\log \weightw \in \Phi_1 \left((u, v)\right)$.
We denote by $P_z$ the functional corresponding to the evaluation of the convolution with
the Poisson kernel at $z \in \mathbb D$;
thus $P_z \psi = \varphi (z)$ for any suitable harmonic function $\varphi$ on $\mathbb D$
having boundary values $\psi$ on $\mathbb T$.
Since functions $\log |u_n|$ are subharmonic,
\begin {equation}
\label {subharmest}
\log |u_n (z, \omega)| \leqslant P_z \log |u_n (\cdot, \omega)| \leqslant P_z [\log \weightw_n (\cdot, \omega)]
\end {equation}
for all $z \in \mathbb D$ and almost all $\omega \in \Omega$.
Passing to the limit in~\eqref {subharmest}
yields $\log |u (z, \omega)| \leqslant P_z [\log \weightw (\cdot, \omega)]$,
which implies that $\weightw \geqslant |u|$ almost everywhere. 
Similarly we obtain $\weightw \geqslant |v|$ and $\weightw \geqslant f \vee g$; thus
$$
\weightw \geqslant |u| \vee |v| \vee f \vee g.
$$
Since $\weightw_n \in 4 B_{X + \lambda Y}$,
by~\cite [Proposition~9] {rutsky2010en} there exists a sequence 
$a_n$ of finite convex combinations of $\{\weightw_j\}_{j \geqslant n}$ such that $a_n \to a \in 4 B_{X + \lambda Y}$
almost everywhere.
On the other hand, weak compactness of $V$ implies
that there exists a sequence $b_n$ of finite convex combinations of $\{\log \weightw_j\}_{j \geqslant n}$
such that $b_n \to \log \weightw$ in $\lclassg {1}$ and hence $b_n \to \log \weightw$ almost everywhere.
A straightforward modification of this construction (see, e.~g., the proof of \cite [Proposition~13] {rutsky2010en})
allows us to assume that the convex combinations $a_n$ and $b_n$ have the same coefficients $\alpha_j^{(n)}$.
By the convexity of the exponential we have
\begin {equation}
\label {expconvex}
\exp b_n = \exp \left(\sum_{j \geqslant n} \alpha_j^{(n)} \log \weightw_j\right) \leqslant
\sum_{j \geqslant n} \alpha_j^{(n)} \weightw_j = a_n.
\end {equation}
Passing to the limit in~\eqref {expconvex} yields $\weightw \leqslant a$ almost everywhere,
and so indeed $\weightw \in 4 B_{X + \lambda Y}$.
Thus the graph of $\Phi_1$ is closed, which implies that $\Phi_1$ is upper semicontinuous and
its values are compact (as closed subsets of $V$); since they are also convex,
we see that $\Phi_1$ is an acyclic map.

Now we endow $(X + \lambda Y)_A$ with the topology of uniform convergence on compact sets in $\mathbb D \times \Omega$
and define a (single-valued) map
$$
\Phi_2 : V \to 4 B_{(X + \lambda Y)_A}
$$
by
\begin {multline}
\label {phi2def}
\Phi_2 (\log \weightw) (z, \omega) =
\exp \left(P_z [\log \weightw (\cdot, \omega)] + i P_z H [\log \weightw (\cdot, \omega)] \right) =
\\
\exp \left(P_z [\log \weightw (\cdot, \omega)] + i Q_z [\log \weightw (\cdot, \omega)] \right)
\end {multline}
for all $\log \weightw \in V$, $z \in \mathbb D$ and a.~e. $\omega \in \Omega$,
where $Q_z$ is the convolution functional with the corresponding conjugate Poisson kernel $Q_z = P_z H$. 
It is easy to see that this map is continuous and $|\Phi_2 (\log \weightw)| = \weightw$ almost everywhere.

Finally, we define a map $\Phi_3 : 4 B_{(X + \lambda Y)_A} \to 2^{B_{X_A} \times \lambda B_{Y_A}}$ by
\begin {equation}
\label {phi3def}
\Phi_3 (h) = \left\{\frac 1 {4 c} (u, v) \mid
(u, v) \in 4 c \left(B_{X_A} \times \lambda B_{Y_A}\right),
h = u + v
\right\}
\end {equation}
for $h \in 4 B_{(X + \lambda Y)_A}$, where $c$ is a constant of the assumed (strong) $\AK$-stability of the couple
$(X, Y)$.
It is easy to see that this map takes nonempty convex values and its graph is closed, so $\Phi_3$ is also an acyclic map.

Now we define the composition map
$F = \Phi_3 \circ \Phi_2 \circ \Phi_1$, which belongs to the class
$\mathfrak B \left(B_{X_A} \times \lambda B_{Y_A}, B_{X_A} \times \lambda B_{Y_A}\right)$.
$F$ is closed and compact as a composition of compact and upper semicontinuous maps.
Thus by Theorem~\ref {parkfptc} there exist some
$u \in B_{X_A}$ and $v \in \lambda B_{Y_A}$ such that $(u, v) \in F ((u, v))$, which means that
there exists some $\weightw \in 4 B_{X + \lambda Y}$ such that $\weightw \geqslant |u| \vee |v| \vee f \vee g$
and $\frac 1 {4 c} \Phi_2 (\log \weightw) = u + v$.
Since $\Phi_2 (\log \weightw)$ is an outer function satisfying 
$\left|\Phi_2 (\log \weightw)\right| = \weightw \geqslant |u| \vee |v| \vee f \vee g$ almost everywhere,
analytic functions
$U = \frac u {\frac 1 {4 c} \Phi_2 (\log \weightw)}$ and
$\frac v {\frac 1 {4 c} \Phi_2 (\log \weightw)} = 1 - U$ are uniformly bounded by $4 c$.
By construction
$$
\|U g\|_X \leqslant 4 c \left\|u\right\|_X \leqslant 4 c = 4 c \|f\|_X
$$
and
$$
\|(1 - U) f\|_Y \leqslant 4 c \|v\|_Y \leqslant 4 c \lambda = 4 c \|g\|_Y.
$$
This concludes the proof of Theorem~\ref {akbakeq}.

\section {Proof of Theorem \ref {akbakeqpc}}

\label {akbakeqpcs}

For any two quasi-normed lattices $X$ and $Y$ on the same measurable space
the set of pointwise products 
$$
X Y = \{ f g \mid f \in X,\break g \in Y\}
$$
is a quasi-normed lattice with the quasi-norm defined by
$$
\|h\|_{X Y} = \inf_{h = f g} \|f\|_X \|g\|_Y.
$$
If both lattices $X$ and $Y$ satisfy the Fatou property then the lattice $X Y$ also has the Fatou property.
If either of the lattices $X$ and $Y$ has order continuous quasi-norm then the quasi-norm of $X Y$ is also order continuous.
By the celebrated Lozanovsky factorization theorem \cite {lozanovsky1969}
$\lclassg {1} = X X'$ for any lattice $X$ satisfying the Fatou property.

For any $\delta > 0$ and a quasi-normed lattice $X$ the lattice $X^\delta$ consists of all measurable functions $f$
such that $|f|^{1/\delta} \in X$ with a quasi-norm $\|f\|_{X^\delta} = \| |f|^{1/\delta} \|_X^\delta$.
For example, $\lclassg {p}^\delta = \lclassg {\frac p \delta}$.
It is easy to see that $(X Y)^\delta = X^\delta Y^\delta$ for any $X$, $Y$ and $\delta$,
and $X^\delta$ naturally inherits many properties from $X$.
For any $0 < \delta \leqslant 1$, if $X$ is a Banach lattice then
$X^\delta$ is also a Banach lattice.
If both $X$ and $Y$ are Banach lattices then for any $0 < \delta < 1$ lattice
$X^{1 - \delta} Y^\delta$, sometimes called the \emph {Calder\'on-Lozanovsky product} of $X$ and $Y$, is also a Banach lattice;
moreover, there is a very useful relation
$(X^{1 - \delta} Y^\delta)' = (X')^{1 - \delta} (Y')^\delta$ (see \cite {calderon1964}, \cite {lozanovsky1969}).
In~\cite {calderon1964} (see also \cite [Chapter~4, Theorem~1.14] {kps})
it was established that Calderón products often describe the complex interpolation spaces between Banach lattices:
$(X_0, X_1)_\theta = X_0^{1 - \theta} X_1^\theta$ 
 if $X_0^{1 - \theta} X_1^\theta$ has order continuous norm (see, e.~g., \cite [Chapter~4, Theorem~1.14] {kps}).

We say that a function $\weightw$ belongs to the \emph {Muckenhoupt class} $\apclass {p}$, $1 < p < \infty$,
with a constant $C$, if
$\weightw \geqslant 0$ a.~e. and
$$
\esssup_{\omega \in \Omega}
\| M \|_{\lclass {p} {\weightw^{-\frac 1 p} (\cdot, \omega)} \to \lclass {p} {\weightw^{-\frac 1 p} (\cdot, \omega)}} \leqslant C;
$$
here $M$ denotes the Hardy-Littlewood maximal operator acting in the first variable.
A lattice $X$ is called \emph {$\apclass {p}$-regular} with constants $(C, m)$ if
for any nonzero $f \in X$ there exists a majorant $\weightu \geqslant |f|$ such that
$\|\weightu\|_X \leqslant m \|f\|_X$ and $\weightu \in \apclass {p}$ with constant $C$.
The well-known relationship between $\apclass {p}$ and $\BMO$ easily implies that
$X$ is a $\BMO$-regular lattice if and only if $X^\delta$ is $\apclass {p}$-regular for some $\delta > 0$.
The notion of $\apclass {p}$-regularity is a very useful refinement of $\BMO$-regularity,
since the sets of the corresponding $\apclass {p}$-majorants are convex;
for more detail see~\cite {rutsky2011en}.

The following proposition generalizes \cite [Proposition~20] {rutsky2013t};
the latter corresponds to the particular case with $r = s = 2$, $Z$ being $2$-convex and both lattices $Z$ and $Z'$
having order continuous norm.
\begin {proposition}
\label {genridiv}
Let $Z$ be a quasi-Banach lattice on $\mathbb T \times \Omega$ having the Fatou property.
Suppose that lattice $V = \left(\lclassg {r}, Z\right)_{\theta, s}$ is $\BMO$-regular
for some $0 < \theta < 1$ and $0 < r, s, < \infty$.
Then $Z$ is also $\BMO$-regular.
\end {proposition}
In order to prove Proposition~\ref {genridiv}, we will reduce it to \cite [Proposition~20] {rutsky2013t}.
The reduction uses fairly straightforward techniques
that appeared in the proof of the latter result; for clarity, we spell out the details.
Indeed, it is well known (see, e.~g., \cite [Theorem~3.2.1] {rolewicz1985}) that a quasi-Banach lattice $Z$
is $q$-convex with some $q > 0$.
By taking $\alpha > 0$ small enough and
replacing the lattice $V$ with $V^\alpha = \left(\lclassg {\frac r \alpha}, Z^\alpha\right)_{\theta, \frac s \alpha}$
(see \cite [Proposition~14] {rutsky2013t} for this equality), $Z$ by $Z^\alpha$, $r$ by $\frac r \alpha$ and
$s$ by $\frac s \alpha$, we may assume that $V$ is $\apclass {2}$-regular,
$Z$ is $2$-convex and $r, s > 2$.
By the reiteration theorem (see \cite [Theorem~3.10.5, Theorem~4.7.2] {bergh}) we have
\begin {equation}
\label {vreit}
V = \left(\lclassg {r}, Z\right)_{\theta, s} =
\left(\lclassg {r}, \left(\lclassg {r}, Z\right)_\beta\right)_{\zeta, s} =
\left(\lclassg {r}, \lclassg {r}^{1 - \beta} Z^\beta\right)_{\zeta, s}
\end {equation}
for all $0 < \beta, \zeta < 1$ such that $\theta = \beta \zeta$.
Such a value of $\beta$ exists
for any $\theta < \zeta < 1$.
Thus we may replace the lattice $Z$ with $\lclassg {r}^{1 - \beta} Z^\beta$ and also assume that
lattices $Z$ and $Z'$ have order continuous norm and $\apclass {2}$-regularity of $V$ holds true for all values of
$\theta$ sufficiently close to~$1$; the corresponding reduction for $\BMO$-regularity follows again
by \cite [Proposition~1.5] {rutsky2011en}.  Furthermore, by the reiteration theorem and
\cite [Proposition~17] {rutsky2013t} lattices
$$
\left(\lclassg {r}, Z\right)_{\theta, p} =
\left(\left(\lclassg {r}, Z\right)_{\theta_0, s}, \left(\lclassg {r}, Z\right)_{\theta_1, s}\right)_{\gamma, p}
$$
are also $\apclass {2}$-regular; here $0 < \theta_0 < \theta_1 < 1$ are sufficiently close to~$1$, $0 < \gamma < 1$,
$\theta = (1 - \gamma) \theta_0 + \gamma \theta_1$, and $1 \leqslant p \leqslant \infty$ is arbitrary.
This implies that we may replace $s$ by any value $1 < s < \infty$, and in particular
we may assume that $s = 2$.
Finally, by applying \cite [Proposition~17] {rutsky2013t} again we see that the lattice
\begin {equation}
\label {threelreit}
\left(\lclassg {t}, \left(\lclassg {r}, Z\right)_{\theta, 2}\right)_{\eta, 2} =
\left(\left(\lclassg {t}, \lclassg {r}\right)_{\theta_2, 2}, Z\right)_{\eta_2, 2} = \left(\lclassg {2}, Z\right)_{\eta_2, 2}
\end {equation}
is $\apclass {2}$-regular; here we take some $1 < t < 2$, find the unique $0 < \theta_2 < 1$ such that
$\frac 1 2 = \frac {1 - \theta_2} t + \frac {\theta_2} r$, and also find the unique $0 < \eta, \eta_2 < 1$
such that $\eta_2 = \theta \eta$ and $(1 - \eta_2) \theta_2 = \eta (1 - \theta)$.  It is easy to see that
these conditions are
satisfied with $\eta_2 = \frac {\theta_2} {\frac 1 \theta - 1 + \theta_2}$ and $\eta = \frac {\eta_2} \theta$.
The first equality in \eqref {threelreit} follows from \cite [Proposition~18] {rutsky2013t}.
This completes the reduction of Proposition~\ref {genridiv} to \cite [Proposition~20] {rutsky2013t}.

Now we are ready to prove Theorem~\ref {akbakeqpc}.  
Since $Y$ is assumed to be $p$-convex,
$Y^p$ is a Banach lattice.  By\footnote {The cited lemma was stated for the strong $\AK$-stability,
but it also works for the usual $\AK$-stability with obvious modifications.}
\cite [Lemma~4] {kisliakov2003}
$\AK$-stability of the couple $(X, Y)$ implies that the couple
$$
\left(X \left[(Y^p)'\right]^{\frac 1 p}, Y \left[(Y^p)'\right]^{\frac 1 p}\right) =
\left(X \left[(Y^p)'\right]^{\frac 1 p}, \lclassg {p}\right)
$$
is also $\AK$-stable; the equality follows from the Lozanovsky factorization formula.
Thus by replacing $X$ with $X \left[(Y^p)'\right]^{\frac 1 p}$ we may assume that
$Y = \lclassg {p}$, which we will do from now on; it is routine to verify that
$\BMO$-regularity of $X \left[(Y^p)'\right]^{\frac 1 p} = \left[X^p (Y^p)'\right]^{\frac 1 p}$
implies that of the couple $(X, Y)$ (see, e.~g., the proof of \cite [Theorem~5.8] {rutsky2011en}).
Now let $Z = \left(X, \lclassg {p}\right)_{\theta, p}$ for some $0 < \theta < 1$ to be determined in a moment.
The Holmstedt formula easily implies that the couple $(Z, \lclassg {p})$ is also $\AK$-stable;
see~\cite [Lemma~1.1] {kisliakov1999}.  The quasi-Banach lattice $X$ is $q$-convex with some $q > 0$.
It is well known (see,\footnote {Observe also that we can easily derive such results for quasi-Banach lattices
from the usual Banach lattice case by exponentiation; see~\cite [Proposition~14] {rutsky2013t}.}
e.~g., \cite [Remark~3 after Proposition~2.g.22] {cbs})
that $Z$ is at least $(r - \varepsilon)$-convex with $\frac 1 r = \frac {1 - \theta} q + \frac \theta p$
and any $\varepsilon > 0$.  Choosing $\theta$ close enough to $1$ thus ensures that $Z$ is a Banach lattice.
We may assume that $p < 2$.  By \cite [Lemma~7] {kisliakov2002en} it follows that the couple
$\left(Z', \lclassg {p'}\right)$ is also $\AK$-stable.
Now we repeat this construction for $\left(Z', \lclassg {p'}\right)$:
taking $Z_1 = \left(Z', \lclassg {p'}\right)_{\theta_1, p'}$ for $0 < \theta_1 < 1$
sufficiently close to $1$ ensures that $Z_1$ is a $2$-convex lattice,
and the couple $\left(Z_1, \lclassg {p'}\right)$
is $\AK$-stable, with $\lclassg {p'}$ also being a $2$-convex lattice.
Thus we may apply \cite [Theorem~2] {rutsky2013t}, which shows that $Z_1$ is $\BMO$-regular.
By Proposition~\ref {genridiv} it follows that lattice $Z'$ is $\BMO$-regular.
Therefore by \cite [Theorem~1.5] {rutsky2011en} lattice $Z$ is also $\BMO$-regular,
and yet another application of Proposition~\ref {genridiv} yields the required $\BMO$-regularity of the lattice $X$.
The proof of Theorem~\ref {akbakeqpc} is complete.

\section {Proof of Theorem~\ref {bmogrilp}}

\label {bmogrilps}

Suppose that under the assumptions of Theorem~\ref {bmogrilp} formula~\eqref {rigip} holds true
for some $0 < \theta < 1$.
Let $X_\alpha = \left(X, \lclassg {q}\right)_{\alpha, p}$ and $Y_\alpha = \left(X_A, \hclassg {q}\right)_{\alpha, p}$
for $0 < \alpha < 1$.  The nontrivial inclusion from~\eqref {rigip} can then be written as
$\left[X_\theta\right]_A \subset Y_\theta$.  By \cite [Theorem~4.7.2] {bergh} we have
$X_\theta = \left(X_{\theta_0}, X_{\theta_1}\right)_{\eta}$ and
$Y_\theta = \left(Y_{\theta_0}, Y_{\theta_1}\right)_{\eta}$ for all $0 < \theta_0 < \theta < \theta_1 < 1$ and
$0 < \eta < 1$ such that $\theta = (1 - \eta) \theta_0 + \eta \theta_1$.
Thus~\eqref {rigip} implies the ``good complex interpolation'' formula
\begin {equation}
\label {gcif}
\left[\left(X_{\theta_0}, X_{\theta_1}\right)_{\eta}\right]_A = \left(Y_{\theta_0}, Y_{\theta_1}\right)_{\eta}
\end {equation}
The assumptions of Theorem~\ref {bmogrilp} imply that lattices $X_\alpha$ are $s$-convex and $t$-concave with some
$1 < s, t < \infty$.  Therefore we may apply \cite [Theorem~5.12] {kalton1994} to \eqref {gcif}, which yields
$\BMO$-regularity of the couple
\begin {equation}
\label {gcbmo1}
\left(X_{\theta_0}, X_{\theta_1}\right) =
\left(\left(X, \lclassg {q}\right)_{\theta_0, p}, \left(X, \lclassg {q}\right)_{\theta_1, p}\right).
\end {equation}
Passing to the dual in \eqref {gcbmo1} if necessary by \cite [Theorem~5.8] {rutsky2011en} allows us to
assume that $q \leqslant 2$.  Raising the couple~\eqref {gcbmo1} to the power $\frac q 2$
by \cite [Proposition~14] {rutsky2013t} and replacing $X$ with $X^{\frac q 2}$ and $p$ with $\frac {2 p} q$
allows us to assume that the couple
$\left(\left(X, \lclassg {2}\right)_{\theta_0, p}, \left(X, \lclassg {2}\right)_{\theta_1, p}\right)$
is $\BMO$-regular, and hence it is $\AK$-stable.  By \cite [Lemma~1.1] {kisliakov1999} and the reiteration theorem
this implies that the couple
\begin {equation}
\label {gcbmo2}
\left(\left(X, \lclassg {2}\right)_{\alpha_0, 2}, \left(X, \lclassg {2}\right)_{\alpha_1, 2}\right)
\end {equation}
is also $\AK$-stable for all $\theta_0 < \alpha_0 < \alpha_1 < \theta_1$, and hence it is $\BMO$-regular
by Theorem~\ref {akbakeqpc}.
By a well-known relation (see, e.~g., \cite [Proposition~11] {rutsky2013t}) we have
\begin {equation}
\label {selfdualformula}
\lclassg {2} = \left(X, X'\right)_{\frac 1 2, 2};
\end {equation}
plugging this formula into \eqref {gcbmo2} and applying the reiteration theorem allows us to rewrite \eqref {gcbmo2} as
\begin {equation*}
\label {gcbmo3}
\left(\left(X, X'\right)_{\beta_0, 2}, \left(X, X'\right)_{\beta_1, 2}\right)
\end {equation*}
with some $0 < \beta_0 < \beta_1 < 1$.
This by~\cite [Theorem~5.8] {rutsky2011en} means that the lattice
\begin {multline*}
\left(X, X'\right)_{\beta_0, 2} \left(X, X'\right)'_{\beta_1, 2} =
\left(X, X'\right)_{\beta_0, 2} \left(X', X''\right)_{\beta_1, 2} =
\\
\left(X, X'\right)_{\beta_0, 2} \left(X', X\right)_{\beta_1, 2} =
\left(X, X'\right)_{\beta_0, 2} \left(X, X'\right)_{1 - \beta_1, 2}
\end {multline*}
is $\BMO$-regular (the order continuity of the norm
assumptions imply that $\left(X, X'\right)'_{\beta_1, 2} = \left(X, X'\right)^*_{\beta_1, 2} =
\left(X^*, {X'}^*\right)_{\beta_1, 2} = \left(X', X''\right)_{\beta_1, 2}$,
and the Fatou property is equivalent to the order reflexivity $X'' = X$),
and so is the lattice
\begin {multline}
\label {gcbmo4}
\left[\left(X, X'\right)_{\beta_0, 2} \left(X, X'\right)_{1 - \beta_1, 2}\right]^{\frac 1 2} =
\\
\left[\left(X, X'\right)_{\beta_0, 2}\right]^{\frac 1 2} \left[\left(X, X'\right)_{1 - \beta_1, 2}\right]^{\frac 1 2} =
\left(X, X'\right)_{\gamma, 2}
\end {multline}
with $\gamma = \frac 1 2 \beta_0 + \frac 1 2 (1 - \beta_1)$;
the last equality follows by \cite [Theorem~4.7.2] {bergh}.
Since $\theta_0$ and $\theta_1$ may take arbitrary values in certain intervals,
we can easily make sure that $\gamma \neq \frac 1 2$.
By passing to the duals in~\eqref {gcbmo4} with the help of \cite [Theorem~1.5] {rutsky2011en} if neccessary we may further
assume that $0 < \gamma < \frac 1 2$, and by the reiteration theorem and \eqref {selfdualformula}
we have
$$
\left(X, X'\right)_{\gamma, 2} = \left(X, \lclassg {2}\right)_{\zeta, 2}
$$
with some $0 < \zeta < 1$.
Finally, Proposition~\ref {genridiv} shows that $X$ is $\BMO$-regular.  The proof of Theorem~\ref {bmogrilp} is complete.

\section {Concluding remarks}

\label {concrem}

The topic of the relationship between $\AK$-stability and $\BMO$-regula\-rity belongs to a large number
of rather general problems,
one of which can be informally stated as follows:
to what extent the interpolation properties of the general couples $(X_A, Y_A)$
are similar to
those of the couples of weighted Lebesgue spaces
$\left(\hclass {\infty} {\weightw_0}, \hclass {\infty} {\weightw_1}\right)$?
As we mentioned in Section~\ref {introduction}, our knowledge of the latter appears to be quite satisfactory.
We can summarize the relevant results from \cite {kisliakov1999}.\footnote {It is also worth mentioning
that except for the partial retractibility (Condition 3) the same is known to be true for the category of
quasi-Banach spaces and all $0 < p_0, p_1 \leqslant \infty$.}
\begin {theorem}
\label {goodinterphp}
Let $\weightw_0$ and $\weightw_1$ be a couple of weights on $\mathbb T \times \Omega$ such that
$\log \weightw_j (\cdot, \omega) \in \lclassg {1}$ for $j \in \{1, 2\}$ and a.~e. $\omega \in \Omega$,
and let $1 \leqslant p_0, p_1 \leqslant \infty$.  The following conditions are equivalent.
\begin {enumerate}
\item
$\left[\mathcal F \left(\left(\lclass {p_0} {\weightw_0}, \lclass {p_1} {\weightw_1} \right)\right)\right]_A
\subset \mathcal F \left(\left(\hclass {p_0} {\weightw_0}, \hclass {p_1} {\weightw_1} \right)\right)$
for all in\-ter\-polation functors $\mathcal F$ in the category of Banach spaces.
\item
$\hclass {p} {\weightw_0^{1 - \theta} \weightw_1^\theta} \subset \left(\hclass {p_0} {\weightw_0}, \hclass {p_1} {\weightw_1} \right)_{\theta, \infty}$,
$\frac 1 p = \frac {1 - \theta} {p_0} +  \frac \theta {p_1}$ for some (equivalently, for all) $0 < \theta < 1$.
\item
$\left(\hclass {p_0} {\weightw_0}, \hclass {p_1} {\weightw_1} \right)$ is a partial retraction of
$\left(\lclass {p_0} {\weightw_0}, \lclass {p_1} {\weightw_1} \right)$.
\item
$\left(\lclass {p_0} {\weightw_0}, \lclass {p_1} {\weightw_1} \right)$ is $\AK$-stable.
\item
$\left(\lclass {p_0} {\weightw_0}, \lclass {p_1} {\weightw_1} \right)$ is $\BMO$-regular.
\item
$\esssup_{\omega \in \Omega} \left\|\log \frac {\weightw_0 (\cdot, \omega)} {\weightw_1 (\cdot, \omega)}\right\|_\BMO < \infty$.
\end {enumerate}
\end {theorem}
If we use this result as a guide, it becomes clear that for a general couple of Banach lattices $(X, Y)$
our knowledge is still somewhat incomplete.
In particular, it is unknown whether $\BMO$-regularity of $(X, Y)$ implies that
$(X_A, Y_A)$ is a partial retraction of~$(X, Y)$.
On the other hand, it is unclear whether a sufficiently general version of $2 \Rightarrow 5$
holds true for the real interpolation; Theorem~\ref {bmogrilp} is a step in this direction.

We also mention that it would be desirable
to revisit the seminal result of N.~Kalton \cite {kalton1994} about the equivalence between
the stability of complex interpolation and $\BMO$-regularity, for clarity and potential generalizations.

\subsection* {Acknowledgement}

The author is grateful to S.~V.~Kisliakov who provided useful remarks to early versions of this paper.

\bibliographystyle {plain}

\bibliography {bmora}

\end {document}